\documentclass[10pt,reqno]{amsart}
\usepackage{amsmath, amssymb}
\usepackage{amsthm}
\usepackage{bm}
\usepackage{tikz}

\newtheorem{theorem}{Theorem}

\theoremstyle{definition}

\newtheorem{example}{Example}
\newcommand{\rats}{\mathbb{Q}}
\newcommand{\Z}{\mathbb{Z}}
\newcommand{\R}{\mathbb{R}}
\newcommand{\Q}{\mathbb{Q}}
\newcommand{\sign}{\operatorname{sgn}}
\newcommand{\lexp}{\operatorname{lexp}}
\newcommand{\supp}{\operatorname{supp}}
\newcommand{\lt}{\operatorname{lt}}

\begin{document}
\bibliographystyle{plain}

\title[Positive Part Extraction for Lattice Paths in the Octant]{Automated Positive Part
Extraction for Lattice Path Generating Functions in the Octant}

\author{Rika Yatchak}
\thanks{Supported by the Austria FWF grant F50-04}
\address{Institute for Algebra, Johannes Kepler University} 
\email{rika.yatchak AT jku.at}

\begin{abstract}
    The question of classifying the nature of the generating functions of restricted lattice
    walks has enjoyed much attention in past years. 
    We prove that a certain class of octant walks have a D-finite
    generating function using the theory of multivariate formal Laurent
    series.
\end{abstract}
\maketitle
\section{Introduction}
We focus our
attention in this paper to the \emph{positive octant}, i.e., the three-dimensional 
integer lattice, restricted so that all coordinates are nonnegative. 

Let $\mathcal{L}=\left(\Z_{\geq 0}\right)^3$. Fix a \emph{model} (also called
\emph{stepset})
$\mathcal{S}\subset\left\{-1,0,1\right\}^3\setminus\{(0,0,0)\}$. A length $m$
$\mathcal{S}$-walk is any walk which starts at $(0,0,0)$ and takes $m$ steps from
$\mathcal{S}$. Let $f(i,j,k;n)$ denote the number of $\mathcal{S}$-walks of
length $n$  ending at point
$(i,j,k)\in\mathcal{L}$. We form the generating
function 
\[F(x,y,z;t)=\sum_{i,j,k,n\geq 0 }f(i,j,k;n)x^iy^jz^kt^n\] 
and ask whether
$F(x,y,z;t)$ is D-finite over $\rats(x,y,z,t)$ in each variable. In other words,
does there exist a nontrivial linear differential equation in each variable
with coefficients in $\rats[x,y,z,t]$? 

There are $2^{3^3}-1= 67{,}108{,}864$ models in the octant. After
choosing a canonical representative for models that are in bijection to each
other for simple reasons and removing all cases that are equivalent to
lower-dimensional problems, we are still left with $10{,}908{,}263$ models
\cite{bacher16}. In \cite{bostan14, bacher16}, these remaining models were
sorted according to properties that sometimes help establish whether the
corresponding generating function is D-finite. One of these properties is the
order of a certain group of rational transformations associated to the
model. If the \emph{group of the model} is finite, then it may be possible to 
establish the D-finite nature of the generating function via the orbit sum
argument of \cite{bousquet10} (summarized in the next section). There are
altogether $2{,}430$ models with a finite group \cite{bacher16,du15}. For $108$ of
them, it was shown in \cite{bostan14} that their corresponding generating
functions are D-finite. We show in the present paper that the orbit sum
argument applies to $1{,}964$ additional models, and that they are therefore also
D-finite. The remaining 358 models are the \emph{zero orbit sum cases}, for
which different methods are needed.

\section{The Orbit Sum Method} 
% a great deal of this section will probably be removed. 
For a fixed model $\mathcal{S}$, we define the \emph{stepset (Laurent)
polynomial} to be 
$P_\mathcal{S}=\sum_{(u,v,w)\in\mathcal{S}} x^u y^v z ^w$. The rational
transformation $\phi_x:\rats(x,y,z)\rightarrow \rats(x,y,z)$ given by
$\phi_x(x,y,z):=
\left(x^{-1}\frac{\sum_{(-1,v,w)\in\mathcal{S}}y^vz^w}{\sum_{(1,v,w)\in\mathcal{S}}y^vz^w},y,z
\right)$
fixes $P_\mathcal{S}$ as well as $y$ and $z$, and is idempotent. The transformations 
$\phi_y$ and $\phi_z$ are defined similarly. Then, the group $G_{\mathcal{S}}$ generated by 
$\phi_x,\phi_y,\phi_z$ via composition  is called the
\emph{group of the model} \cite{bousquet10,bostan14}.
        For $g\in G_{\mathcal{S}}$, $\sign(g) = 1$ if $g$
        is a composition of an even number of generators, and
        $-1$ otherwise. 

        Most models in the octant have an infinite group, but for the 2{,}430 cases in the
octant with a finite group, we can attempt to
use the \emph{orbit sum method} to prove that their
corresponding generating functions are D-finite. The details of this technique for
the octant appear in \cite{bostan14}; we reproduce a brief sketch here.

The starting point of the technique is the \emph{functional equation}, which
can be obtained using an inclusion-exclusion argument. The left-hand side of this
equation only involves $F(x,y,z;t)$, but the right-hand side involves
\emph{sections} or \emph{specializations} of $F(x,y,z;t)$ such as $F(x,y,0;t)$ and 
$F(x,0,0;t)$ \cite[Eqn. 7]{bostan14}. Using the
group of the model, we form the \emph{orbit sum}, which allows us to eliminate
the sections on the right-hand side and leaves us with: 
\begin{equation}
        \sum_{g\in G}
        \sign(g)g(xyz)(F(g(x),g(y),g(z);t))=\frac{1}{1-tP_\mathcal{S}}\sum_{g\in G}
\sign(g)g(xyz)
\label{eq:os}
\end{equation}

The generating function $F(x,y,z;t)$ is an element of $\Q[x,y,z][[t]]$. Thus,
$F_g:=F(g(x),g(y),g(z);t)\in \Q(x,y,z)[[t]]$ is a power series in $t$ with
coefficients $a_i\in \Q(x,y,z)$. We regard each of these $a_i$ as a
\emph{multivariate formal Laurent series} rather than a rational function (see
Section \ref{sec:MFLS}). 
Then, the positive part extraction $[x^{>0}y^{>0}z^{>0}]$ is 
well-defined, and we apply it to obtain an expression with only $F(x,y,z;t)$ on
the right-hand side

\begin{equation}
        xyzF(x,y,z;t)=[x^{>0}y^{>0}z^{>0} ]\frac{1}{1-tP_{\mathcal{S}}(x,y,z)}\sum_{g\in G}
        \sign(g)g(xyz)\label{eq:pp} 
\end{equation}
The above equation holds if every monomial in 
$F_g$ has at least one negative
component in its exponent for each non-identity element $g\in G_{\mathcal{S}}$.  
Then, $F(x,y,z;t)$ can be written
as the diagonal of a rational series, and is therefore
D-finite \cite{bostan17,lipshitz88}. 

The orbit sum technique gives a uniform method for proving that some octant
models have D-finite generating functions. However, it is not without its
limitations. If there are one or more non-identity $g\in G_{S}$ such that
$F(g(x),g(y),g(z);t))$
contributes to the positive part on the left-hand side, we do not obtain an
expression for $xyzF(x,y,z;t)$ alone, and therefore cannot conclude that
$F(x,y,z;t)$ is D-finite. For such models, a different proof technique is
needed. Included among these models are the so-called \emph{zero orbit sum}
cases, for which the right-hand side of Equation \ref{eq:os} is $0$. In the
quadrant, these models correspond exactly to the cases where the generating
function is algebraic \cite{bousquet10,bostan10}. In the octant, there are two
subclasses among the zero orbit sum cases: \emph{Hadamard} and
``mysterious". Bostan et al provide an alternative proof technique for
octant Hadamard cases in \cite{bostan14}: all of the zero orbit sum Hadamard
models they consider have D-finite generating functions. The remaining 170 models are those 
which were deemed ``mysterious" in \cite{bacher16}. These
cases have yet to be definitively classified, but there is strong
computational evidence that at least some of them have non-D-finite generating
functions \cite{bacher16}.

In this paper, we identify models for which $xyzF(x,y,z;t)$ is the only
surviving term on the left-hand side of Equation \ref{eq:os} after the positive part extraction,
and prove automatically that their generating functions are D-finite via the
orbit sum argument. 

\section{Multivariate formal Laurent series\label{sec:MFLS}}
We recall here the bare essentials of the theory of multivariate formal Laurent series, as
formalized by Aparicio-Monforte and Kauers \cite{aparicio12}. In order to
ensure that multiplication of series is well-defined, we consider 
multivariate series with support contained in a \emph{line-free} cone.

Let $k$ be a field, $x_1,\dots ,x_n$ be $n$ indeterminates, and $C\subset \R^n$ be a
line-free cone. $k_C[[x_1,\dots,x_n]]:=\left\{f(\bm{x})=\sum_{\bm{k}}
        a_{\bm{k}}\bm{x}^{\bm{k}} \mid \supp
f(\bm{x})\subseteq C\right\}$ with the usual addition and Cauchy product
multiplication is an integral domain. 
If, on the other hand, we fix a specific additive order on $\Z^n$, we can
form the sets $k_{\preceq}[[\bm{x}]]:=\bigcup_{C\in \mathcal{C}}
k_C[[\bm{x}]]$ and
$k_{\preceq}((\bm{x})):=\bigcup_{\bm{v}\in\Z^p}\bm{x}^{\bm{v}}k_\preceq
[[\bm{x}]]$, where $\bm{x}=(x_1,\dots ,x_p)$ with $x_i$ indeterminate for
every~$i$, 
and $\mathcal{C}$ is the set of all cones $C\subset \R^p$
\emph{compatible} with $\preceq$. The condition that $\preceq$ be additive
ensures that $k_\preceq [[\bm{x}]]$ is a ring and that
$k_\preceq((\bm{x}))$ is a field \cite[Thm. 15]{aparicio12}. We call
$k_\preceq[[\bm{x}]]$ a multivariate formal Laurent series ring and $k_{\preceq}((\bm{x}))$ the
field of the multivariate formal Laurent series ring. 
We define the \emph{leading exponent} of an element $f(\bm{x})$ to
be 
$\lexp_{\preceq} f(\bm{x})=\min_{\preceq}( \supp f(\bm{x}))\in \Z^{p}$, and
the \emph{leading term} $\lt_\preceq f(\bm{x})=\bm{x}^{\lexp_\preceq
f(\bm{x})}$.

The following theorem is essential for the automatic positive part extraction
we will describe in the next section: 
\begin{theorem}[{\cite[Thm. 17]{aparicio12}}] 
        Let $C\subset R^{q}$ be a line-free cone and $f(\bm{y})\in
        k_C[[y]]$. Let $\preceq$ be an additive order on $\Z^p$ and
        $a_1(\bm{x}),\dots ,a_p(\bm{x})\in k_\preceq((x))\setminus\{0\}$. Let $M\in
        \Z^{p\times q}$ be the matrix whose $i^\text{\tiny th}$ column consists
        of the leading exponent $\lexp(a_i(\bm{x}))$. Let $C'\subset R$
        be a cone containing $MC:=\left\{M\bm{x}\mid \bm{x}\in C\right\}\subset
        \R^p$ as well as
        $\supp\left(a(\bm{x})/\lt(a_i(\bm{X}))\right)$ for every
        $i$ $\in \{1,\dots ,q\}$. Suppose that $C\cap \ker M=\left\{\bm{0}\right\}$ and that
        $C'$ is line-free. Then, $f\left(a_1(\bm{x}),\dots,a_p(\bm{x})\right)$
        is well-defined and belongs to the ring $k_{C'}[[\bm{x}]]$. \label{thm:cone}
\end{theorem}

After the positive part extraction is
applied with the help of this theorem, the next step is to find a differential
equation for $F(x,y,z;t)$ and possibly an expression for $F(x,y,z;t)$ itself.
The theory of multivariate formal Laurent series is also useful for this step:  
in \cite{bostan17}, for example, this theory
is used to prove computationally guessed annihilating differential
operators for the sections $F(x,0;t)$ and $F(0,y;t)$ of certain quadrant models. 
These operators lead to an annihilating differential operator for $F(x,y;t)$,
as well as explicit expressions for $F(x,y;t)$ in 
terms of hypergeometric functions. 

\section{Application to Lattice Paths\label{sec:app}}
% The formatting is klugey here but it works
For a given $\mathcal{S}$, we apply Theorem \ref{thm:cone} $|G_\mathcal{S}|$ times. The
support of $F(x,y,z;t)$ can be shown to be contained in $C = \langle
(1,0,0,1),(0,1,0,1),(0,0,1,1),(0,0,0,1),\\(1,1,0,1),(1,0,1,1),(0,1,1,1),(1,1,1,1)\rangle$. 
Fix a non-identity $g\in G_{\mathcal{S}}$. We set $f(\bm{y})=F(x,y,z;t)$, $a_1(\bm{x}):=g(x)$, $a_2(\bm{x}):=g(y)$,
$a_3(\bm{x})=g(z)$, and $a_4(\bm{x})=t$. 
To create an additive order on $\Z^p$, we collect all polynomials $q_i$ that appear
in the numerator or denominator of the $a_i$. For each of these polynomials
$q_i$ we choose a leading term, and check that this choice of leading term is
compatible with the cone. If so, we check the conditions of 
Theorem \ref{thm:cone} are fulfilled. If they are, we conclude that the composition
$F_g=F(g(x),g(y),g(z);t)$ is valid, and also obtain a cone $C'_g$ with the property that $F_g\in
k_{C'_g}[[\bm{x}]]$. Next, define $B$ to be the smallest cone containing $C$ and $C'_g$. We check that $B$
is line-free, and that 
$\supp(g(xyz)F_g)\cap C=\emptyset$. Then, we have proven that
$[x^{>0}y^{>0}z^{>0}]g(xyz)F_g=0$.  

If we experience a failure at any of the above steps, we simply choose different leading terms for
the $q_i$ and try again. The properties of $G_\mathcal{S}$ ensure that there are not many polynomials
$q_i$ for which we can choose leading terms, and that each $q_i$ can only
contain a small number of monomials. Thus, it is computationally feasible to
check every combination of leading terms that is compatible with $C$.

We repeat this process for every non-identity $g\in G$, setting $B$ to be equal to the
smallest cone containing the previous cone $B$ and the cone $C'_g$. If the
process terminates successfully, we obtain a cone $B$ with the property that
$g(xyz)F_g\in k_\preceq((\bm{x}))$ for every $g$. We also know that 
the only contribution to the left-hand side of the orbit sum with support intersecting 
$C$ is the element $xyzF(x,y,z)$.  
That is, the positive part extraction step yields an expression for $xyzF(x,y,z;t)$ 
alone that allows us to conclude that
$F(x,y,z;t)$ is D-finite. 

\begin{example}
    Let $\mathcal{S}=\{(-1,-1,0),(-1,0,1),(-1,1,-1),
(0,-1,1),(0,0,-1),\\(0,1,0),(1,0,0)\}$.
$G_\mathcal{S}=\langle\phi_x,\phi_y,\phi_z\rangle\cong D_{12}$, with
$\phi_x=\left(\frac{yz^2+y^2+z}{xyz},y,z\right),\phi_y=\left(x,\frac{z}{y},
z\right),\phi_z=\left(x,y,\frac{y}{z}\right)$. For this example, the only
polynomial that appears as a term in a rational transformation is 
$p=yz^2+y^2+z$. Choose $\lexp_\preceq(p)=(0,1,2,0)$. For each $g\in
G$ we now attempt to apply Theorem \ref{thm:cone}.  
Consider for example the element $\phi_x\in G_{\mathcal{S}}$.
$M=\left(\begin{array}{cccc}-1 & 0 & 0 & 0\\
                            0 & 1 & 0 & 0\\
                            1 & 0 & 1 & 0\\
                            0 & 0 & 0 & 1 \end{array}\right)$.

We have that $C\cap\text{ker}(M)=\{\bm{0}\}$. Additionally, 
$C'_{\phi_x} = \langle (0,-2,1,0),\\(0,-1,2,0),(0,0,0,1)\rangle$. Thus,
$F_{\phi_x}=\frac{yz^2+y^2+z}{xyz}yzF\left(\frac{yz^2+y^2+z}{xyz},y,z;t\right)$ is well-defined, and additionally
$\supp(\phi_x(xyz)F_{\phi_x})\cap C = \emptyset$. Moreover, the smallest cone
$B$ containing both
$C$ and $C'_{\phi_x}$ is line-free. Continuing this process for each $g\in
G_{\mathcal{S}}$, we obtain the cone $C_G=\langle( 1,  1,  1, 1),( 1, -1, 2, 1),
( 0,  1, -2, 0),(-1, -1,  3, 1),\\( 0, -1, -1, 0),
(-1,  1,  2, 1) \rangle$, which is line-free. The left-hand side of the
orbit sum equation is well-defined in $k_{C_G}[[\bm{x}]]$, and $F(x,y,z;t)$
is the only term that survives the operation $\supp\left(g(xyz)\cdot
F(g(x),g(y),g(z);t)\right)\cap C$ for $g\in G_\mathcal{S}$, and $F(x,y,z;t)$ is
D-finite.
\end{example}

\begin{example} 
    Next we consider a case for which the right-hand side of Equation
    \ref{eq:os} is equal to zero. Let $\mathcal{S}=\{(-1, -1, -1), (-1, 0, 0), (-1, 0, 1), (-1, 1, 0),
                (1, -1, 0),\\
     (1, 0, -1), (1, 1, 1)\}$. $G_\mathcal{S}\cong D_{12}$ as before, but is
     generated by different rational transformations $\phi_x,\phi_y,\phi_z$
     than in the previous example. In particular,
     we have 
     $\phi_z\phi_y\phi_z = (x,z,y)$. Since $\lt_\preceq(x)=x$, $\lt_\preceq
     (y)=y$ and $\lt_{\preceq} (z)=z$ for any choice of $\preceq$, we have
     $M_{\phi_z\phi_y\phi_z}=\left(\begin{array}{cccc}
                                            1 & 0 & 0 & 0\\
                                            0 & 0 & 1 & 0\\
                                            0 & 1 & 0 & 0\\
                                            0 & 0 & 0 & 1 \end{array}\right)$
    Since $C'_{\phi_z\phi_y\phi_z}=C$, the term $xyzF(x,z,y)$ on the left-hand
    side of Equation \ref{eq:os} survives the
    positive part extraction, and we cannot proceed.
\end{example}

\section{Octant cases with finite group}
Using the technique outlined in the previous section, we obtain the
following theorem. 

\begin{theorem} All $2{,}072$ three-dimensional octant models with finite
        group and nonzero orbit sum are D-finite. \label{thm:dfin}
\end{theorem}

We say that a model is three-dimensional if the condition that $\sum_{s\in\mathcal{S}}a_s 
s \geq (0,0,0)$ is a
truly three-constraint problem. A model defined in the octant need not
be three-dimensional; see \cite[\S 2.1,\S 7.1]{bostan14} for a more extended
discussion.

Theorem \ref{thm:dfin} includes the 108 cases that were already
proven in \cite{bostan14}. In \cite{bacher16}, it was conjectured that the
$1{,}964$ models 
covered in Theorem \ref{thm:dfin} could be proved D-finite by the orbit 
sum method. The contribution of the current
work is to confirm this conjecture. The proof is fully automatic: our
implementation of the process outlined in Section \ref{sec:app} uses Sage
\cite{sage} and requires only the model $\mathcal{S}$ as input. 

The remaining 358 models not covered by Theorem \ref{thm:dfin} are exactly 
those for which Equation \ref{eq:os} has a vanishing right-hand side. This
shows that the only cases for which the orbit sum technique fails at the
positive part extraction step are those that are equivalent to
two-dimensional cases with multiplicities, which are discussed in
\cite{bostan14, kauers15c}. 

%\section{Conclusion} 
%We have shown rigorously and automatically that a large class of octant
%cases with finite associated group have D-finite generating functions. One 
%direction for further work is to classify models for which  
%orbit sum is zero.  In the quarter plane, these cases are in direct
%correspondence with those cases where the generating function is not only
%D-finite, but algebraic. In the octant, it appears that things are not so
%simple: the guessed asymptotics in \cite{bacher16} provide strong
%computational evidence that some of the models with an associated orbit sum
%of $0$ are non-D-finite, but no rigorous proof of this result is yet known.

\section*{Acknowledgements} 
I would like to thank my adviser Manuel Kauers for many helpful discussions
on this topic. 
\bibliography{refs}

\begin{thebibliography}{10}

\bibitem{aparicio12}
Ainhoa {Aparicio Monforte} and Manuel Kauers.
\newblock Formal {L}aurent series in several variables.
\newblock {\em Expositiones Mathematicae}, 31(4):350--367, 2013.

\bibitem{bacher16}
Axel Bacher, Manuel Kauers, and Rika Yatchak.
\newblock Continued classification of 3{D} lattice walks in the positive
  octant.
\newblock {\em ArXiv e-prints}, November 2015.
\newblock To appear in proceedings of FPSAC 2016.

\bibitem{bostan14}
Alin Bostan, Mireille Bousquet-M{\'e}lou, Manuel Kauers, and Stephen Melczer.
\newblock On 3-dimensional lattice walks confined to the positive octant.
\newblock {\em Annals of Combinatorics}, 2014.
\newblock to appear.

\bibitem{bostan17}
Alin Bostan, Fr{\'e}d{\'e}ric Chyzak, Mark van Hoeij, Manuel Kauers, and Lucien
  Pech.
\newblock Hypergeometric expressions for generating functions of walks with
  small steps in the quarter plane.
\newblock {\em European Journal of Combinatorics}, 61:242 -- 275, 2017.

\bibitem{bostan10}
Alin Bostan and Manuel Kauers.
\newblock The complete generating function for {G}essel walks is algebraic.
\newblock {\em Proceedings of the AMS}, 138(9):3063--3078, 2010.
\newblock with an appendix by Mark van Hoeij.

\bibitem{bousquet10}
Mireille Bousquet-M{\'e}lou and Marni Mishna.
\newblock Walks with small steps in the quarter plane.
\newblock {\em Contemporary Mathematics}, 520:1--40, 2010.

\bibitem{sage}
The~Sage Developers.
\newblock {\em {S}age {M}athematics {S}oftware ({V}ersion 6.8)}, 2015.
\newblock {\tt http://www.sagemath.org}.

\bibitem{du15}
Daniel~K. Du, Qing-Hu Hou, and Rong-Hua Wang.
\newblock Infinite orders and non-{D}-finite property of 3-dimensional lattice
  walks.
\newblock {\em Electron. J. Combin.}, 23(3):Paper 3.38, 15, 2016.

\bibitem{kauers15c}
Manuel Kauers and Rika Yatchak.
\newblock Walks in the quarter plane with multiple steps.
\newblock In {\em Proceedings of FPSAC'15}, DMTCS, pages 35--36, 2015.

\bibitem{lipshitz88}
Leonard Lipshitz.
\newblock The diagonal of a {$D$}-finite power series is {$D$}-finite.
\newblock {\em J. Algebra}, 113(2):373--378, 1988.

\end{thebibliography}
\end{document}